# ON THE CLASSIFICATION OF CERTAIN PIECEWISE LINEAR AND DIFFERENTIABLE MANIFOLDS IN DIMENSION EIGHT AND AUTOMORPHISMS OF $\#_{i=1}^{b}(S^2 \times S^5)^\dagger$

ALEXANDER SCHMITT

ABSTRACT. In this paper, we will be concerned with the explicit classification of closed, oriented, simply-connected spin manifolds in dimension eight with vanishing cohomology in the odd dimensions. The study of such manifolds was begun by Stefan Müller. In order to understand the structure of these manifolds, we will analyze their minimal handle presentations and describe explicitly to what extent these handle presentations are determined by the cohomology ring and the characteristic classes. It turns out that the cohomology ring and the characteristic classes do not suffice to reconstruct a manifold of the above type completely. In fact, the group $\mathrm{Aut}_0\bigl(\#_{i=1}^{b}(S^2 \times S^5)\bigr)/\mathrm{Aut}_0\bigl(\#_{i=1}^{b}(S^2 \times D^6)\bigr)$ of automorphisms of $\#_{i=1}^{b}(S^2 \times S^5)$ which induce the identity on cohomology modulo those which extend to $\#_{i=1}^{b}(S^2 \times D^6)$ acts on the set of oriented homeomorphy classes of manifolds with fixed cohomology ring and characteristic classes, and we will be also concerned with describing this group and some facts about the above action.

## 1. INTRODUCTION

The classification of topological manifolds up to homeomorphy is an extremely interesting and important problem. Let us restrict to the case of closed (i.e., compact without boundary) and oriented simply connected manifolds. As a general classification scheme, surgery theory [1] solves this problem for manifolds within a given homotopy type, e.g., the one of a sphere. Another approach to the classification "up to finite indeterminacy", using rational homotopy theory, is due to Sullivan [34]. Nevertheless, there are only a few explicit results which characterize the oriented homeomorphy type of a manifold in terms of easily computable invariants. They usually require a lot of simplifying assumptions such as high connectedness [36]. In this paper, we will consider *even cohomology manifolds* (or *E-manifolds*, for short) in dimension eight by which we mean closed, oriented, simply connected, piecewise linear or smooth manifolds all of whose odd dimensional homology groups with integer coefficients vanish. The universal coefficient theorem implies that all homology groups of an E-manifold are without torsion. Moreover, since $H^3(X, \mathbb{Z}_2) = 0$ for an E-manifold, two E-manifolds of dimension at least 6 are homeomorphic (as topological manifolds) if and only if they are piecewise linearly homeomorphic [16].

Though the class of E-manifolds is fairly restricted, it still contains many interesting examples from various areas of mathematics, such as the piecewise linear manifolds underlying the toric manifolds from Algebraic Geometry [4] or the polygon spaces [12], to name a few. So far, E-manifolds have been classified up to dimension 6. Of course, in dimension 2 there is only $S^2$, in dimension 4, there is the famous classification result of Freedman various interesting aspects of which are discussed in [15], and finally in dimension 6, the classification was achieved by Wall [37] and Jupp [14]. Various applications of the latter result to Algebraic Geometry are surveyed in [26]. Finally, we refer to [2], [3], and [29] for the determination of projective algebraic structures on certain 6- and 8-dimensional E-manifolds.

## ACKNOWLEDGMENT

My thanks go to J.-C. Hausmann for his interest in this work and for pointing out to me the fact that an E-manifold of dimension eight is actually not determined by its classical invariants (which was asserted by me in the first version of this paper).

## 2. STATEMENT OF THE RESULT

We now discuss the main result of this note, namely the classification of E-manifolds of dimension 8 with vanishing second Stiefel-Whitney class in the form of an exact sequence of pointed sets. This result was motivated by the work [24]. In order to state it in a more elegant form, we will work with based manifolds. By a *based piecewise linear (smooth) E-manifold*, we mean a triple $(X,\underline{x},\underline{y})$, consisting of a piecewise linear (smooth) E-manifold $X$ and bases $\underline{x} = (x_1,...,x_{b_2(M)})$ for $H^2(X,\mathbb{Z})$ and $\underline{y} = (y_1,...,y_{b_4(M)})$ for $H^4(X,\mathbb{Z})$. Recall that E-manifolds are by definition oriented, so that the above data specify a basis for $H^*(X,\mathbb{Z})$, such that the bases for $H^i(X,\mathbb{Z})$ and $H^{8-i}(X,\mathbb{Z})$ are dual to each other w.r.t. the cup product. An *isomorphism between piecewise linear (smooth) based E-manifolds* $(X,\underline{x},\underline{y})$ and $(X',\underline{x}',\underline{y}')$ is an orientation preserving piecewise linear (smooth) isomorphism $f\colon X \longrightarrow X'$ with $f^*(\underline{x}') = \underline{x}$ and $f^*(\underline{y}') = \underline{y}$. Denote by $\mathfrak{I}^{\mathrm{PL}(\mathscr{C}^\infty)}(b,b')$ the set of isomorphy classes of piecewise linear (smooth) based E-manifolds $(X,\underline{x},\underline{y})$ of dimension eight with vanishing second Stiefel-Whitney class, $b_2(X) = b$, and $b_4(X) = b'$.

2.1. **The classical invariants.** In the terminology of [24], the classical invariants of an E-manifold consist of its cohomology ring, the Stiefel-Whitney classes, the Wu classes, the Pontrjagin classes, the Euler class, the Steenrod squares, the reduced Steenrod powers, and the Pontrjagin powers. For an eight dimensional E-manifold $X$ with vanishing second Stiefel-Whitney class, the main result of [24] states that the classical invariants are fully determined by the following invariants:

1. The cup product map
$$\begin{aligned} \delta_X \colon S^2 H^2(X,\mathbb{Z}) &\longrightarrow H^4(X,\mathbb{Z}) \\ x \otimes x' &\longmapsto x \cup x'. \end{aligned}$$

2. The intersection form
$$\begin{aligned} \gamma_X \colon S^2 H^4(X,\mathbb{Z}) &\longrightarrow \mathbb{Z} \\ y \otimes y' &\longmapsto (y \cup y')[X]. \end{aligned}$$

Here, $[X] \in H_8(X,\mathbb{Z})$ is the fundamental class determined by the orientation.
3. The first Pontrjagin class $p_1(X) \in H^4(X,\mathbb{Q})$.

*Remark* 2.1. The above invariants are not independent. By associativity of the cohomology ring, the following relation holds

$$\delta_X^*(\gamma_X) \in S^4 H^2(X,\mathbb{Z})^\vee, \tag{1}$$

i.e.,

$$\gamma_X\big(\delta_X(x_1 \otimes x_2) \otimes \delta_X(x_3 \otimes x_4)\big) = \gamma_X\big(\delta_X(x_1 \otimes x_3) \otimes \delta_X(x_2 \otimes x_4)\big),$$

for all $x_1, x_2, x_3, x_4 \in H^2(X, \mathbb{Z})$.

Furthermore, one has

**Proposition** ([24], Prop. A.7 or Cor. 3.14 below). *For every element $y \in H^4(X, \mathbb{Z})$*

$$p_1(X) y \equiv 2y^2 \mod 4. \tag{2}$$

Observe that this implies $p_1(X) \in H^4(X, \mathbb{Z})$.

If $X$ is in addition differentiable, one has for every integral lift $W \in H^2(X, \mathbb{Z})$ of $w_2(X)$

$$3p_1(X)^2 - 14 p_1(X) W^2 + 7 W^4 \equiv 12 \operatorname{Sign}(\gamma_X) \mod 2688. \tag{3}$$

Müller has also shown that these relations imply all other relations among the classical invariants of $X$ [24]. Conversely, a piecewise linear manifold $X$ the invariants of which obey relation (3) admits a differentiable structure [18], [24].

We are led to the following algebraic concept: A *system of invariants of type* $(b, b')$ is a triple $(\delta, \gamma, p)$, consisting of

- a homomorphism $\delta \colon S^2 \mathbb{Z}^{\oplus b} \longrightarrow \mathbb{Z}^{\oplus b'}$,
- a unimodular symmetric bilinear form $\gamma \colon S^2 \mathbb{Z}^{\oplus b'} \longrightarrow \mathbb{Z}$, and
- an element $p \in \mathbb{Z}^{\oplus b'}$.

Denote by $Z(b, b')$ the set of systems of invariants of type $(b, b')$.

Now, let $(X, \underline{x}, \underline{y})$ be a based eight dimensional E-manifold. This defines a set of invariants $Z_{(X,\underline{x},\underline{y})} := (\delta_X, \gamma_X, p_1(X))$ of type $(b_2(X), b_4(X))$. Thus, we have natural maps

$$Z^{\mathrm{PL}(\mathscr{C}^\infty)}(b, b') \colon \mathfrak{J}^{\mathrm{PL}(\mathscr{C}^\infty)}(b, b') \longrightarrow Z(b, b')$$
$$[X, \underline{x}, \underline{y}] \longmapsto Z_{(X,\underline{x},\underline{y})}.$$

It will be the concern of our paper to understand the maps $Z^{\mathrm{PL}(\mathscr{C}^\infty)}$ as well as possible. The first result can be easily derived from Wall's work [36] and deals with the case $b = 0$. It will be proved in detail in Section 4.1.

**Theorem 2.2.** i) *The map $Z^{\mathrm{PL}}(0, b')$ is injective. Its image consists precisely of those elements which satisfy the relations (1) and (2).*

ii) *Given two smooth based E-manifolds $(X, \underline{y})$ and $(X', \underline{y}')$ with $b_2(X) = 0 = b_2(X')$ and $Z_{(X,\underline{y})} = Z_{(X',\underline{y}')}$, there exists an exotic 8-sphere $\Sigma$ such that $(X \# \Sigma, \underline{y})$ and $(X', \underline{y}')$ are smoothly isomorphic. In particular, the fibres of $Z^{\mathscr{C}^\infty}$ have cardinality at most two. The image of $Z^{\mathscr{C}^\infty}$ consists exactly of those elements which satisfy the relations (1), (2), and (3).*

2.2. **Manifolds with trivial cup form $\delta_X$.** In addition to describing the explicit geometric meaning of the system of invariants of an E-manifold $X$ with $w_2(X) = 0$, we will describe those manifolds $X$ for which the cup form $\delta_X$ vanishes.

For any $b > 0$, let $\mathscr{FC}_b^{\mathrm{PL}(\mathscr{C}^\infty)}$ be the group of isotopy classes of piecewise linear (smooth) embeddings of $b$ disjoint copies of $S^5 \times D^3$ into $S^8$. The following result will be established in Section 3.7.

**Proposition 2.3.**
$$\mathrm{FL}_b \quad := \quad \mathscr{FC}_b^{\mathscr{C}^\infty} \quad \cong \quad \mathscr{FC}_b^{\mathrm{PL}}$$

Given an element $[l] \in \mathrm{FL}_b$, we can perform surgery along the link $l$ and get a smooth based E-manifold $(X(l), \underline{x}(l))$ with $w_2(X) = 0$ and $b_4(X) = 0$ which is well defined up to smooth isomorphy of based manifolds.

We will also use the following notation: Fix a pair $(\gamma, p) \in Z(0, b')$ which satisfies relation (2) (and (3)) and denote by $\mathfrak{J}^{\mathrm{PL}(\mathscr{C}^\infty)}(b, \gamma, p)$ the set of isomorphy classes of based piecewise linear (smooth) E-manifolds $(X, \underline{x}, \underline{y})$ with $w_2(X) = 0$, $b_2(X) = b$, $\gamma_X = \gamma$, and $p_1(X) = p$. Pick a three-connected piecewise linear (smooth) based E-manifold $(X_0, \underline{y}_0)$ with $\gamma_{X_0} = \gamma$ and $p_1(X_0) = p$. In the smooth case, let $\vartheta^8 \cong \mathbb{Z}_2$ [17] be the group of exotic 8-spheres, and set $\vartheta(X_0) := \vartheta^8$, if $X_0$ is not diffeomorphic to $X_0 \# \Sigma$, $\Sigma$ a generator for $\vartheta^8$, and $\vartheta(X_0) := \{[S^8]\} \subset \vartheta^8$ otherwise. Now, we define maps

$$K^{\mathrm{PL}}(b, \gamma, p) \colon \mathrm{FL}_b \quad \longrightarrow \quad \mathfrak{J}^{\mathrm{PL}}(b, \gamma, p)$$
$$[l] \quad \longmapsto \quad \left[ X(l) \# X_0, \underline{x}(l), \underline{y}_0 \right]$$

and

$$K_{X_0}^{\mathscr{C}^\infty}(b, \gamma, p) \colon \mathrm{FL}_b \oplus \vartheta(X_0) \quad \longrightarrow \quad \mathfrak{J}^{\mathscr{C}^\infty}(b, \gamma, p)$$
$$([l], [\Sigma]) \quad \longmapsto \quad \left[ X(l) \# X_0 \# \Sigma, \underline{x}(l), \underline{y}_0 \right].$$

In $\mathfrak{J}^{\mathrm{PL}(\mathscr{C}^\infty)}(b, \gamma, p)$, we mark the class $\left[ (\#_{i=1}^b (S^2 \times S^6)) \# X_0, \underline{x}, \underline{y}_0 \right]$ where $\underline{x}$ comes from the natural basis of $H^2(\#_{i=1}^b (S^2 \times S^6), \mathbb{Z})$. Then, our main result is the following

**Theorem 2.4.** i) *For every $b > 0$ and every pair $(\gamma, p)$ which satisfies the relation (2),*

$$\begin{array}{ccccccc} \{1\} & \longrightarrow & \mathrm{FL}_b & \stackrel{K^{\mathrm{PL}}(b,\gamma,p)}{\longrightarrow} & \mathfrak{J}^{\mathrm{PL}}(b, \gamma, \varphi) & \longrightarrow & \mathrm{Hom}(S^2 \mathbb{Z}^b, Z^{b'}) \\ 1 & \longmapsto & [\text{trivial link}] & & [X, \underline{x}, \underline{y}] & \longmapsto & \delta_X \end{array}$$

*is an exact sequence of pointed sets.*

ii) *For every $b > 0$ and every pair $(\gamma, p)$ which satisfies the relations (2) and (3),*

$$\begin{array}{ccccccc} \{1\} & \longrightarrow & \mathrm{FL}_b \oplus \vartheta(X_0) & \stackrel{K_{X_0}^{\mathscr{C}^\infty}(b,\gamma,p)}{\longrightarrow} & \mathfrak{J}^{\mathscr{C}^\infty}(b, \gamma, \varphi) & \longrightarrow & \mathrm{Hom}(S^2 \mathbb{Z}^b, Z^{b'}) \\ 1 & \longmapsto & [\text{trivial link}] & & [X, \underline{x}, \underline{y}] & \longmapsto & \delta_X \end{array}$$

*is an exact sequence of pointed sets.*

The proof will be given in Section 4.2 and 4.3.

*Remark* 2.5. i) In the PL setting, we will show that the inclusion of $\text{FL}_b$ into $\mathfrak{J}^{\text{PL}}(b,\gamma,\varphi)$ extends to an action of $\text{FL}_b$ on $\mathfrak{J}^{\text{PL}}(b,\gamma,\varphi)$.

ii) On all the sets occuring in Theorem 2.4, there are natural $(\text{GL}_b(\mathbb{Z}) \times \text{GL}_{b'}(\mathbb{Z}))$-actions, and the maps are equivariant for these actions. Therefore, by forming the equivalence classes w.r.t. these actions, we get the classification of E-manifolds with vanishing second Stiefel-Whitney class up to orientation preserving piecewise linear (smooth) isomorphy.

iii) We will discuss in Section 3.7 the structure of the group $\text{FL}_b$. It turns out that it is finite if and only if $b = 1$. It follows easily that the set of $\text{GL}_b(\mathbb{Z})$-equivalence classes in $\text{FL}_b$ is infinite for $b \geq 2$. Thus, the cohomology ring and the characteristic classes classify E-manifolds of dimension eight up to finite indeterminacy only if the second Betti number is at most one.

Starting point of our proof of the above results will be the theory of minimal handle decompositions of Smale which states that $X$ can be obtained from $D^8$ by first attaching $b_2(X)$ 2-handles, then $b_4(X)$ 4-handles, then $b_2(X)$ 6-handles and finally one 8-handle. At each step, the attachment will be determined by the isotopy class of a certain framed link in a 7-manifold, and we will first explain how to read off the isotopy class of the attaching links for the 2- and 4-handles from the invariants.

## 3. PRELIMINARIES

We collect in this paragraph the background material and some preliminary results which we will use in our proof. Most of the results are by now standard results from various parts of algebraic, differential, and piecewise linear topology.

### 3.1. **The structure of manifolds: handle attachment and surgery.**

Let $M$ be an $m$-dimensional manifold with boundary. Suppose we are given an embedding $f\colon S^{\lambda-1} \times D^{m-\lambda} \longrightarrow \partial M$. We then define

$$M' \quad := \quad M \cup_f \left(D^\lambda \times D^{m-\lambda}\right)$$

and say that $M'$ is obtained from $M$ by the *attachment of a $\lambda$-handle along $f$*. Moreover, $f(S^{\lambda-1} \times \{0\})$ is called the *attaching sphere*, $D^\lambda \times \{0\}$ the *core disc*, and $\{0\} \times S^{m-\lambda-1}$ the *belt sphere*. We will often simply write $M' = M \cup H^\lambda$.

*Remark* 3.1. i) If $M$ is assumed to be differentiable and $f$ to be a differentiable embedding, handle attachment can be described in such a way that the resulting manifold is again differentiable (see [17], VI, §§6&8), i.e., no "smoothing of the corners" is required.

ii) If $M$ is oriented, then $M'$ will inherit an orientation which is compatible with the given orientation of $M'$ and the natural orientation of $D^\lambda \times D^{m-\lambda}$, if and only if $f$ *reverses* the orientations.

The next operation we consider was introduced by Milnor [21] and Wallace [38] and goes back to Thom. For this, let $N$ be a manifold of dimension $m-1$ and $f\colon S^{\lambda-1} \times D^{m-\lambda} \longrightarrow N$ an embedding. Denote by $\overline{f}$ the restriction of $f$ to $S^{\lambda-1} \times S^{m-\lambda-1}$, and set

$$\chi(N,f) \quad := \quad \left(N \setminus \text{int} f(S^{\lambda-1} \times D^{m-\lambda})\right) \cup_{\overline{f}} \left(D^\lambda \times S^{m-\lambda-1}\right).$$

We say that $\chi(N,f)$ is constructed from $N$ by *surgery along $f$*. Informally speaking, we remove from $N$ a $(\lambda-1)$-sphere with trivial normal bundle and replace it with an $(m-\lambda-1)$-sphere, again with trivial normal bundle.

*Remark* 3.2. i) If $N$ is oriented, then $f$ has to be orientation preserving for $\chi(N,f)$ to inherit a natural orientation from those of $N$ and $D^\lambda \times S^{m-\lambda-1}$. This is because $S^{\lambda-1} \times S^{m-\lambda-1}$ inherits the reversed orientation as boundary of $N \setminus \text{int } f(S^{\lambda-1} \times D^{m-\lambda})$.

ii) The operations of handle attachment and surgery are closely related: Let $M$ be an $m$-dimensional manifold with boundary $N := \partial M$ and $f\colon S^{\lambda-1} \times D^{m-\lambda} \longrightarrow N$ an embedding. Now, attach a $\lambda$-handle along $f$ in order to obtain $M'$. Then, $\partial M' = \chi(N,f)$.

We will also perform a "surgery in pairs". For this, $N$ is assumed to be an $(m-1)$-dimensional manifold and $K$ to be a submanifold of dimension $k-1$. Assume that, for some $\lambda \leq k$, we are given an embedding $f\colon S^{\lambda-1} \times D^{m-\lambda} \longrightarrow N$ which induces an embedding $f^* := f_{|S^{\lambda-1} \times D^{k-\lambda}}\colon S^{\lambda-1} \times D^{k-\lambda} \longrightarrow K$. Then, $\chi(K,f^*)$ is naturally contained as a submanifold in $\chi(N,f)$.

The next result is a special case of the "minimal presentation theorem" of Smale [31] and is crucial for the explicit analysis of the structure of a manifold.

**Theorem 3.3.** *Let $X$ be a closed simply connected manifold of dimension $m \geq 6$ with torsion free homology. Then, there exists a sequence of submanifolds*

$$D^m \cong W_0 \subset W_1 \subset W_2 \subset \cdots \subset W_m = X,$$

*such that $W_i$ is obtained from $W_{i-1}$ by attaching $b_i(X)$ $i$-handles, $i = 1,\ldots,m$.*

*Moreover, for any such sequence, there exists a* dual sequence

$$\overline{W}_0 \subset \overline{W}_1 \subset \cdots \subset \overline{W}_m = X,$$

*such that the attaching $(i-1)$-spheres in $\overline{W}_{i-1}$ coincide with the belt spheres in $W_{m-i}$, $i = 1,\ldots,m$.*

*Proof.* For differentiable manifolds, an attractive presentation of the relevant material is contained in Chapter VII and VIII of [17]. In the piecewise linear category, handle decompositions are discussed in [27] (cf. also [13]). However, the statement concerning the number of handles is not explicitly given there. Nevertheless, one verifies that the necessary tools (such as Whitney lemma and handle sliding) are also proved in [27]. □

*Remark* 3.4. i) Retracting all $\lambda$-handles to their core discs, starting with $\lambda = 0$, yields a CW-complex which is homotopy equivalent to $X$ (cf. [27], p. 83).

ii) Observe that, by i), a handle decomposition as in Theorem 3.3 yields a *preferred basis* for $H_*(X,\mathbb{Z})$. By renumbering, orientation reversal in the attaching spheres, and handle sliding, one can obtain any basis of $H_*(X,\mathbb{Z})$ as the preferred basis of a handle decomposition ([17], (1.7), p. 148)

iii) If $X$ comes with an orientation, we may assume that $D^m$ is orientation preservingly embedded and that all attaching maps are orientation reversing.

3.2. **Consequences for E-manifolds of dimension eight with $w_2 = 0$.** Let $X$ be a piecewise linear (smooth) E-manifold of dimension eight with $w_2(X) = 0$. The first observation concerns the structure of $W_2$.

**Lemma 3.5.** *One has $W_2 \cong \#_{i=1}^{b}(S^2 \times D^6)$.*

*Proof.* The manifold $W_2$ is an $(8,1)$-handle body and as such homeomorphic to the boundary connected sum of $b$ $D^6$-bundles over $S^2$ ([17], §11, p. 115). As $\pi_1(SO(4)) \cong \mathbb{Z}_2$ and we have requested $w_2(X) = 0$, the claim follows. □

The next consequence is

*The manifold $W_4$ is determined by a framed link of $b_4(X)$ three dimensional spheres in $\#_{i=1}^{b}(S^2 \times S^5)$.*

We will address the classification of such links below.

The third consequence is

**Lemma 3.6.** i) $\partial W_4 \cong \#_{i=1}^{b}(S^2 \times S^5)$.
  ii) *The manifold $X$ is of the form $W_4 \cup_f \#_{i=1}^{b}(S^2 \times D^6)$ where*

$$f\colon \#_{i=1}^{b}(S^2 \times S^5) \longrightarrow \partial W_4$$

*is a piecewise linear (smooth) isomorphism, such that $f_*$ maps the canonical basis of $H_2(\#_{i=1}^{b}(S^2 \times S^5), \mathbb{Z})$ to the preferred basis of $H_2(\partial W^4, \mathbb{Z})$.*

*Proof.* i) This follows because $\partial W_4 \cong \partial \overline{W}_2$. ii) follows because $X = W_4 \cup \overline{W}_2$, and $\overline{W}_2 \cong \#_{i=1}^{b}(S^2 \times D^6)$, by Lemma 3.5. □

3.3. **Homotopy vs. isotopy.** By Theorem 3.3, the manifold is determined by the ambient isotopy classes of the attaching maps. However, the topological invariants of the manifold give us at best their homotopy classes. It is, therefore, important to have theorems granting that this is enough. In the setting of differentiable manifolds, we have

**Theorem 3.7** (Haefliger [6], [7]). *Let S be a closed differentiable manifold of dimension n and M an m-dimensional differentiable manifold without boundary. Let $f\colon S \longrightarrow M$ be a continuous map and $k \geq 0$, such that*

$$\pi_i(f)\colon \pi_i(S) \longrightarrow \pi_i(M)$$

*is an isomorphism for $0 \leq i \leq k$ and surjective for $i = k+1$. Then, the following is satisfied:*

1. *If $m \geq 2n - k$ and $n > 2k + 2$, then $f$ is homotopic to a differentiable embedding.*
2. *If $m > 2n - k$ and $n \geq 2k + 2$, then two differentiable embeddings of S into M which are homotopic are also ambient isotopic.*

In the setting of piecewise linear manifolds, similar results hold true. We refer to Haefliger's survey article [9]. For our purposes, the result stated below will be sufficient.

**Theorem 3.8.** *Suppose S is a closed n-dimensional manifold, M an m-dimensional manifold without boundary, and $f\colon S \longrightarrow M$ a continuous map. Assume one has*

- $m - n \geq 3$.
- *S is $(2n - m + 1)$-connected.*
- *M is $(2n - m + 2)$-connected.*

*Then:*

1. *$f$ is homotopic to an embedding.*

2. *Two embeddings which are homotopic to f are ambient isotopic.*

*Proof.* The theorem of Irwin ([27], Thm. 7.12 and Ex. 7.14, [13], Thm. 8.1) yields 1. and that $f_1$ and $f_2$ as in 2. are concordant. But, since $m-n \geq 3$, concordance implies ambient isotopy ([13], Chap. IX). □

**Corollary 3.9.** *Let $S := S^3$ and $M$ a simply connected differentiable or piecewise linear manifold of dimension 7 without boundary. Then, $\pi_3(M)$ classifies differentiable and piecewise linear embeddings, respectively, of $S^3$ into $M$ up to ambient isotopy.*

**Corollary 3.10** (Zeeman's unknotting theorem [39])**.** *For every $m,n$ with $m-n \geq 3$, any piecewise linear embedding of $S^n$ into $S^m$ is isotopic to the standard embedding.*

3.4. **Some 4-dimensional CW-complexes.** By Remark 3.4, a handle decomposition of $X$ gives us a CW-complex which is homotopy equivalent to $X$. The following discussion will enable us to understand the 4-skeleton of that complex.

Let $W := S^2 \vee \cdots \vee S^2$ be the $b$-fold wedge product of 2-spheres. Suppose $X$ is the CW-complex obtained by attaching a 4-cell to $W$ via the map $g \in \pi_3(W)$. The Hilton-Milnor theorem ([30], Thm. 7.9.4) asserts

$$\pi_3(W) = \bigoplus_{i=1}^{b} \pi_3(S^2) \oplus \bigoplus_{1 \leq i < j \leq b} \pi_3(S^3).$$

Choosing the standard generators for $\pi_3(S^2)$ and $\pi_3(S^3)$, we can describe $g$ by a tuple $(l_i, i=1,...,b; l_{ij}, 1 \leq i < j \leq b)$ of integers. These integers determine the cohomology ring of $X = W \cup_g D^4$ as follows:

**Proposition 3.11.** *Let $y \in H^4(X, \mathbb{Z})$ be the generator of $H^4(X, \mathbb{Z})$ given by the attached 4-cell and $x_1,...,x_b$ the canonical basis of $H^2(X, \mathbb{Z}) = H^2(W, \mathbb{Z})$. Then*

$$\begin{aligned} x_i \cup x_j &= l_{ij} \cdot y, & 1 \leq i < j \leq b, \\ x_i \cup x_i &= l_i \cdot y, & i=1,...,b. \end{aligned}$$

This is proved like [22], (1.5), p. 103. We recall the proof in the following example.

*Example* 3.12. We treat the case $b=2$. Consider the embedding

$$\iota: S^2 \vee S^2 \hookrightarrow S^2 \times S^2 \hookrightarrow \mathbb{CP}^\infty \times \mathbb{CP}^\infty.$$

The standard basis for $H^4(\mathbb{CP}^\infty \times \mathbb{CP}^\infty, \mathbb{Z}) \cong \mathbb{Z}^{\oplus 3}$ is given by the elements $y_1$, $y_2$, $y_3$ obtained from attaching $D^4$ via $(1,0;0)$, $(0,0;1)$, and $(0,1;0)$, respectively. Let $h: D^4 \longrightarrow D^4 \vee D^4 \vee D^4$ be the canonical map followed by

$$\bigl(\vartheta \cdot x \longmapsto \vartheta \cdot m_{l_1}(x)\bigr) \vee \bigl(\vartheta \cdot x \longmapsto \vartheta \cdot m_{l_{12}}(x)\bigr) \vee \bigl(\vartheta \cdot x \longmapsto \vartheta \cdot m_{l_2}(x)\bigr).$$

Here, $m_k$ stands for a representative of $[k \cdot \mathrm{id}_{S^3}] \in \pi_3(S^3)$ and $D^4 = \{\vartheta \cdot x \,|\, x \in S^3, \vartheta \in [0,1]\}$. Now, $h$ and $\iota$ glue to a map $f: X \longrightarrow \mathbb{CP}^\infty \times \mathbb{CP}^\infty$, and

$$\begin{aligned} f^*: H^4(\mathbb{CP}^\infty \times \mathbb{CP}^\infty, \mathbb{Z}) &\longrightarrow H^4(X, \mathbb{Z}) \\ a_1 y_1 + a_2 y_2 + a_3 y_3 &\longmapsto (a_1 l_1 + a_2 l_{12} + a_3 l_2) y, \end{aligned}$$

so that the assertion follows from the naturality of the cup-product.

### 3.5. Pontrjagin classes and $\pi_3(SO(4))$.

Vector bundles of rank 4 over $S^4$ are classified by elements in $\pi_3(SO(4))$. In our setting, such vector bundles will appear as normal bundles. We recall, therefore, the description of that group and relate it to Pontrjagin classes and self intersection numbers.

First, look at the natural map $\pi_3(SO(4)) \longrightarrow \pi_3(SO(4)/SO(3)) = \pi_3(S^3)$. This map has a splitting ([32], §22.6) which induces an isomorphism

$$\pi_3(SO(4)) = \pi_3(SO(3)) \oplus \pi_3(S^3).$$

Let $\alpha_3$ be the generator for $\pi_3(SO(3)) \cong \mathbb{Z}$ from [32], §22.3, and $\beta_3 := [\mathrm{id}_{S^3}] \in \pi_3(S^3)$, so that we obtain the isomorphism $\mathbb{Z} \oplus \mathbb{Z} \longrightarrow \pi_3(SO(4))$, $(k_1, k_2) \longmapsto k_1\alpha_3 + k_2\beta_3$. Finally, the kernel of the map $\pi_3(SO(4)) \longrightarrow \pi_3(SO)$ to the stable homotopy group is generated by $-\alpha_3 + 2\beta_3$ ([32], §23.6), whence [23], (20.9), implies

**Proposition 3.13.** *Let $E$ be the vector bundle over $S^4$ defined by the element $k_1\alpha_3 + k_2\beta_3 \in \pi_3(SO(4))$. Then*

$$p_1(E) = \pm(2k_1 + 4k_2).$$

**Corollary 3.14.** *Let $f: S^4 \longrightarrow M$ be a differentiable embedding of $S^4$ into the differentiable 8-manifold $M$. Let $E := f^*T_M/T_{S^4}$ be the normal bundle. Then, the self intersection number $s$ of $f(S^4)$ in $M$ satisfies*

$$2s \equiv p_1(E) \mod 4.$$

*Proof.* If $E$ is given by the element $k_1\alpha_3 + k_2\beta_3 \in \pi_3(SO(4))$, then $s = k_2$ ([17], (5.4), p. 72). Since $p_1(E) = \pm(2k_2 + 4k_1)$, the claim follows. □

### 3.6. Links of 3-spheres in $\#_{i=1}^{b}(S^2 \times S^5)$.

If $X$ is a closed E-manifold of dimension 8 with $w_2(X) = 0$, then $W_2 := \#_{i=1}^{b}(S^2 \times D^6)$, $b = b_2(X)$, by Lemma 3.5. Thus, $W_4$ is determined by a framed link of 3-spheres in $\partial W_2 = \#_{i=1}^{b}(S^2 \times S^5)$. Therefore, we will now classify such links.

So, let $W := \#_{i=1}^{b}(S^2 \times S^5)$ be a $b$-fold connected sum. We can choose $b$ disjoint 2-spheres $S_i^2$, $i = 1,...,b$, embedded in $W$ and representing the natural basis of $H_2(W, \mathbb{Z})$. One checks that the homotopy type of $W$ is given up to dimension 4 by the $b$-fold wedge product $S^2 \vee \cdots \vee S^2$. Suppose we are given a link of $b'$ three dimensional spheres, i.e., we are given $b'$ differentiable embeddings $g_i: S^3 \longrightarrow W$, $i = 1,...,b'$, with $g_i(S^3) \cap g_j(S^3) = \emptyset$ for $i \neq j$.

By the transversality theorem ([17], IV.(2.4)), one sees that we may assume $S_i^2 \cap g_j(S^3) = \emptyset$ for all $i$ and $j$.

By Corollary 3.9, the ambient isotopy class of the embedding $g_k$ is determined by the element $\varphi_k := [g_k] \in \pi_3(W_k)$, $W_k := W \setminus \bigcup_{j \neq k} g_j(S^3)$, $k = 1,...,b'$. We clearly have (compare [8])

$$\pi_3(W_k) = \pi_3\bigl(\underbrace{S^2 \vee \cdots \vee S^2}_{b\times} \vee \underbrace{S^3 \vee \cdots \vee S^3}_{(b'-1)\times}\bigr),$$

so that the Hilton-Milnor theorem yields

$$\pi_3(W_k) = \bigoplus_{i=1}^{b} \pi_3(S^2) \oplus \bigoplus_{1 \leq i < j \leq b} \pi_3(S^3) \oplus \bigoplus_{j \neq k} \pi_3(S^3).$$

Hence, we write $\varphi_k$ as a tuple of integers:

$$\varphi_k = \left(l_i^k, i=1,...,b; l_{ij}^k, 1 \leq i < j \leq b; \lambda_{kj}, j \neq k\right).$$

Observe that, for $j \neq k$, $\varphi_k$ is mapped under the natural homomorphism

$$\pi_3(W_k) \longrightarrow H_3(W_k, \mathbb{Z}) \longrightarrow H_3\bigl(W \setminus g_j(S^3), \mathbb{Z}\bigr)(\cong \mathbb{Z})$$

to the image of the fundamental class of $S^3$ under $g_{j*}$. Thus, $\lambda_{kj}$ is just the "usual" linking number of the spheres $g_k(S^3)$ and $g_j(S^3)$ in $W$ (compare [8]).

3.7. **Links of 5-spheres in $S^8$.** Let $\mathscr{FC}_b^{\mathrm{PL}(\mathscr{C}^\infty)}$ be as before, and let $\mathscr{C}_b^{\mathrm{PL}(\mathscr{C}^\infty)}$ be the group of isotopy classes of piecewise linear (smooth) embeddings of $b$ disjoint copies of $S^5$ into $S^8$. For $b=1$, these groups are studied in [10], [19], and [20]. A brief summary with references of results in the case $b > 1$ is contained in Section 2.6 of [11]. We will review some of this material below.

**Proposition 3.15.** *We have $\mathscr{FC}_1^{\mathscr{C}^\infty} \cong \mathscr{FC}_1^{\mathrm{PL}} \cong \mathbb{Z}_2$.*

*Proof.* Since $\pi_5(\mathrm{SO}(3)) \cong \mathbb{Z}_2$, the standard embedding of $S^5$ into $S^8$ with its two possible framings provides an injection of $\mathbb{Z}_2$ into $\mathscr{FC}_1^{\mathrm{PL}(\mathscr{C}^\infty)}$. By Zeeman's unknotting theorem 3.10, the map $\mathbb{Z}_2 \longrightarrow \mathscr{FC}_1^{\mathrm{PL}}$ is an isomorphism. As remarked in Section 2.6 of [11], $\mathscr{FC}_1^{\mathrm{PL}}$ is isomorphic to $\mathscr{F}\vartheta$, the group of h-cobordism classes of framed submanifolds of $S^8$ which are homotopy 5-spheres. Moreover, by [10] and [19], there is an exact sequence

$$\cdots \longrightarrow \vartheta^6 \longrightarrow \mathscr{FC}_1^{\mathscr{C}^\infty} \longrightarrow \mathscr{F}\vartheta \longrightarrow \vartheta^5 \longrightarrow \cdots.$$

As the groups $\vartheta^5$ and $\vartheta^6$ of exotic 5- and 6-spheres are trivial [17], our claim is settled. □

Let $\mathrm{L}_b \subset \mathscr{C}_b^{\mathscr{C}^\infty}$ be the subgroup of those embeddings for which the restriction to each component is isotopic to the standard embedding. As observed in Section 2.6 of [11], Zeeman's unknotting theorem 3.10 implies that $\mathrm{L}_b = \mathscr{C}_b^{\mathrm{PL}}$. The following settles Proposition 2.3:

**Corollary 3.16.**

$$\mathscr{FC}_b^{\mathscr{C}^\infty} \cong \mathscr{FC}_b^{\mathrm{PL}} \cong \mathrm{L}_b \oplus \bigoplus_{i=1}^b \mathbb{Z}_2.$$

For the group $\mathrm{L}_b$, Theorem 1.3 of [11] provides a fairly explicit description as an extension of abelian groups. For this, consider the $b$-fold wedge product $\bigvee_{i=1}^b S^2$ of 2-spheres together with its inclusion $i\colon \bigvee_{i=1}^b S^2 \hookrightarrow \bigtimes_{i=1}^b S^2$ into the $b$-fold product of 2-spheres. Finally, let $p_i\colon \bigvee_{i=1}^b S^2 \longrightarrow S^2$ be the projection onto the $i$-th factor, $i=1,...,b$. Set, for $m=1,2,...$,

$$\Lambda_{b,j}^m := \mathrm{Ker}\bigl(\pi_m(p_j)\colon \pi_m(\bigvee_{i=1}^b S^2) \longrightarrow \pi_m(S^2)\bigr), \quad j=1,...,b,$$

$$\Lambda_b^m := \bigoplus_{j=1}^b \Lambda_{b,j}^m$$

and

$$\Pi_b^m := \operatorname{Ker}\bigl(\pi_m(i)\colon \pi_m(\bigvee_{i=1}^b S^2) \longrightarrow \bigoplus_{i=1}^b \pi_m(S^2)\bigr),$$

and define

$$w_b^m \colon \Lambda_b^m \longrightarrow \Pi_b^{m+1}$$

on $\Lambda_{b,j}^m$ by $w_b^m(\alpha) := [\alpha, \iota_i]$. Here, $[.,.]$ stands for the Whitehead product inside the homotopy groups of $\bigvee_{i=1}^b S^2$ and $\iota_i \colon S^2 \hookrightarrow \bigvee_{i=1}^b S^2$ for the inclusion of the $i$-th factor, $i = 1, ..., b$. Theorem 1.3 of [11] yields in our situation

**Theorem 3.17.** *There is an exact sequence of abelian groups*

$$0 \longrightarrow \operatorname{Coker}(w_b^6) \longrightarrow L_b \longrightarrow \operatorname{Ker}(w_b^5) \longrightarrow 0.$$

We remark that the formulas of Steer [33] might be used for the explicit computation of Whitehead products and thus for the determination of $\operatorname{Coker}(w_b^6)$ and $\operatorname{Ker}(w_b^5)$. The free part of $L_b$, e.g., can be obtained quite easily. We confine ourselves to prove the following important fact.

**Corollary 3.18.** *The group $L_b$ has positive rank for $b \geq 2$.*

*Proof.* Let $\mathbb{L}_b := \bigoplus_{l \geq 1} \mathbb{L}_{b,l}$ be the free graded Lie algebra with $\mathbb{L}_{b,1} := \bigoplus_{i=1}^b \mathbb{Z} \cdot e_i$. For $l = 2, 3, ...$, let $e_1^l, ..., e_{d_l}^l$ be a basis for $\mathbb{L}_{b,l}$ consisting of iterated commutators of the $e_i$. By assigning $\iota_i$ to $e_i$, every iterated commutator $c \in \mathbb{L}_{b,l}$ in the $e_i$ defines an element $\alpha(c) \in \pi_{l+1}(\bigvee_{i=1}^b S^2)$.

To settle our claim, it is certainly sufficient to show that $\operatorname{Coker}(w_b^6)$ has positive rank. Now, by the Hilton-Milnor theorem

$$\Pi_b^7 \cong \bigoplus_{l=3}^7 \bigoplus_{k=1}^{d_{l-1}} \pi_7(S^l) \cdot \alpha(e_k^{l-1}).$$

Note that $\pi_7(S^l)$ is finite for $l \notin \{4, 7\}$ (see [32] and [35] for the explicit description of those groups). The Hopf fibration $S^7 \longrightarrow S^4$ [32], on the other hand, yields a decomposition $\pi_7(S^4) \cong \pi_6(S^3) \oplus \pi_7(S^7) \cong \mathbb{Z}_{12} \oplus \mathbb{Z}$. Therefore, it will suffice to show that the free part of $\Lambda_b^6$ is mapped to $\bigoplus_{j=1}^{d_6} \pi_7(S^7) \cdot \alpha(e_j^6)$. For $j = 1, ..., b$, we have

$$\Lambda_{b,j}^6 \cong \bigoplus_{i \neq j} \pi_6(S^2) \cdot \iota_i \oplus \bigoplus_{l=3}^6 \bigoplus_{k=1}^{d_{l-1}} \pi_6(S^l) \cdot \alpha(e_k^{l-1}).$$

The group $\pi_6(S^l)$ is finite for $l < 6$, and we obviously have $[\alpha(e_k^5), \iota_j] = \alpha([e_k^5, e_j])$. If we expand the commutator $[e_k^5, e_j]$ in the basis $e_1^6, ..., e_{d_6}^6$, we find an expansion for $[\alpha(e_k^6), \iota_j]$ in terms of the $\alpha(e_k^6)$. $\square$

**Corollary 3.19.** *The set of $\operatorname{GL}_b(\mathbb{Z})$-equivalence classes of elements in $L_b$ is infinite for $b \geq 2$.*

*Proof.* We have seen that the $\operatorname{GL}_b(\mathbb{Z})$-set $\mathbb{L}_{b,3}$ is contained in the $\operatorname{GL}_b(\mathbb{Z})$-set $L_b$. The $\operatorname{GL}_b(\mathbb{Z})$-action on $\mathbb{L}_{b,3}$ origins from a homomorphism $\operatorname{GL}_b(\mathbb{Z}) \longrightarrow \operatorname{GL}(\mathbb{L}_{b,3}) := \operatorname{Aut}_{\mathbb{Z}}(\mathbb{L}_{b,3})$. In particular, any matrix $g \in \operatorname{GL}_b(\mathbb{Z})$ preserves the absolute value of the determinant of any $d_3$ elements in $\mathbb{L}_{b,3}$. This implies, for instance, that $a \cdot e_1^3$ and $b \cdot e_1^3$ cannot lie in the same $\mathrm{GL}_b(\mathbb{Z})$-orbit, if $0 \leq a < b$.  □

## 4. Proof of Theorem 2.2 and Theorem 2.4

From now on, $X$ stands for an eight dimensional E-manifold with $w_2(X) = 0$.

4.1. **Proof of Theorem 2.2.** The classification result for 3-connected E-manifolds of dimension eight is a special case of a result of Wall's [36] and can be easily obtained with the methods described in [17], VII, §12. Let us recall the details, because we will need them later on.

We fix a basis $\underline{b}$ for $H_4(X,\mathbb{Z})$ and let $\underline{y}$ be the dual basis of $H^4(X,Z)$. Then, there is a handle presentation $X = D^8 \cup H_1^4 \cup \cdots \cup H_{b'}^4 \cup D^8$ with $\underline{b}$ as the preferred basis. The manifold $T := D^8 \cup H_1^4 \cup \cdots \cup H_{b'}^4$ is determined by the ambient isotopy class of a framed link of 3-spheres in $S^7$, having $b'$ components. Let us first look at such a link, forgetting the framing, i.e., suppose we are given embeddings $g_i \colon S^3 \longrightarrow S^7$ with $S_i \cap S_j = \varnothing$ for $i \neq j$, $S_i := g_i(S^3)$, $i = 1,...,b'$. By 3.9, we may assume that the $g_i$ are differentiable. Observe that the normal bundles of the $S_i$ are trivial. We equip $S_i$ with the orientation induced via $g_i$ by the standard orientation of $S^3$ and the normal bundle of $S_i$ with the orientation which is determined by requiring that the orientation of $S_i$ followed by the one of its normal bundle coincide with the orientation of $S^7$. Therefore, a 3-sphere $F_i$ which bounds the fibre of a tubular neighborhood of $S_i$ in $S^7$ inherits an orientation and thus provides a generator $e_i$ for $H_3(S^7 \setminus S_i, \mathbb{Z}) \cong \mathbb{Z}$, $i = 1,...,b'$. For $i \neq j$, the image of the fundamental class $[S_i]$ in $H_3(S^7 \setminus S_j, \mathbb{Z})$ is of the form $\lambda_{ij} \cdot e_j$. The integer $\lambda_{ij}$ is called *the linking number of $S_i$ and $S_j$*.

For $i = 1,...,b'$, the manifold $S^7 \setminus \bigcup_{j \neq i} S_j$ is up to dimension 5 homotopy equivalent to $\bigvee_{j \neq i} F_j$, and

$$\pi_3\bigl(S^7 \setminus \bigcup_{j \neq i} S_j\bigr) \cong \pi_3\bigl(\bigvee_{j \neq i} F_j\bigr) \cong \bigoplus_{j \neq i} H_3(S^7 \setminus S_j, \mathbb{Z}).$$

Under this identification, we have $[g_i] = \sum_{j \neq i} \lambda_{ij} \cdot e_j$. The $[g_i]$ determine the ambient isotopy class of the given link (3.9), and we deduce

**Proposition 4.1.** *The linking numbers $\lambda_{ij}$, $1 \leq i < j \leq b'$, determine the given link up to ambient isotopy.*

The sphere $S_i$ bounds a 4-dimensional disc $D_i^-$ in $D^8$, $i = 1,...,b'$, which we equip with the induced orientation. We may, furthermore, assume that the $D_i^-$ intersect transversely in the interior of $D^8$. Then, the $\lambda_{ij}$ coincide with the intersection numbers $D_i^- . D_j^-$, $1 \leq i < j \leq b'$. For an intuitive proof (in dimension 4), see [28], p. 67. Now, every disc $D_i^-$ is completed by the core disc $D_i^+$ of the $i$-th 4-handle to an embedded 4-sphere $\Sigma_i$ in $T$, $i = 1,...,b'$, and, since all the core discs are pairwise disjoint, the $\lambda_{ij}$ coincide with the intersection numbers $\Sigma_i . \Sigma_j$, $1 \leq i < j \leq b'$. Finally, $X$ is obtained by gluing an 8-disc to $T$ along $\partial T$, and the spheres $\Sigma_i$ represent the elements of the chosen basis $\underline{b}$ of $H_4(X,\mathbb{Z})$. Identifying the intersection ring with the cohomology ring of $X$ via Poincaré-duality, we see

**Corollary 4.2.** *The linking numbers $\lambda_{ij}$ coincide with the cup products $(y_i \cup y_j)[X]$, $1 \leq i < j \leq b'$, i.e., the link of the attaching spheres is determined up to ambient isotopy by the basis $\underline{b}$ and the cup products.*

As we have remarked before, the normal bundles of the $S_i$ are trivial, whence there exist embeddings $f_i^0 \colon S^3 \times D^4 \longrightarrow S^7$ with $f^0_{i|S^3 \times \{0\}} = g_i$, $i = 1, ..., b'$. From the uniqueness of tubular (in differential topology) or regular (in piecewise linear topology) neighbourhoods, every other embedding $f_i \colon S^3 \times D^4 \longrightarrow S^7$ with $f_{i|S^3 \times \{0\}} = g_i$ is ambient isotopic to one of the form $f_i^{[h_i]} := \big((x,y) \longmapsto (x, h_i \cdot y)\big)$, $[h_i] \in \pi_3(\mathrm{SO}(4))$, $i = 1, ..., b'$. Corollary 3.14 implies that we can choose the $f_i^0$, $i = 1, ...b'$, in such a way that the following holds:

**Lemma 4.3.** *Suppose $T$ is obtained by attaching 4-handles along $f_i^{[h_i]}$ with $[h_i] = k_1^i \alpha_3 + k_2^i \beta_3$, $i = 1, ..., b'$, then*

$$\Sigma_i . \Sigma_i \;=\; k_2^i \quad \text{and} \quad p_1(T_{T|\Sigma_i}) = \pm\big(2k_2^i + 4k_1^i\big).$$

This shows that also the framed link used for constructing $T$ and $X$ is determined by the system of invariants associated to $(X, \underline{y})$, proving the injectivity in Part i) in the theorem. Moreover, the assertion about the fibres in Part ii) is clear.

Conversely, given a system $Z$ of invariants in $Z(0, b')$, satisfying relation (2), there exists a based 3-connected manifold $(X, \underline{y})$ realizing $Z$. Indeed, by the above identification of the invariants, $Z$ determines a framed link in $S^7$ and thus the manifold $T := D^8 \cup H_1^4 \cup \cdots \cup H_{b'}^4$. The boundary of $T$ is a 7-dimensional homotopy sphere ([17], (12.2), p. 119) and, therefore, piecewise linearly homeomorphic to $S^7$. Hence, $X = T \cup_{S^7} D^8$ is a piecewise linear manifold with the desired system of invariants, settling Part i). If, in addition, relation (3) holds, the work [18] grants that $X$ will carry a smooth structure (compare Theorem A.4 of [24]), finishing the proof of Part ii). □

4.2. **The determination of $W_4$ in the general case.** We have a handle decomposition $W_0 \subset W_2 \subset W_4 \subset W_6 \subset X$ of $X$ providing preferred bases $\underline{b}$ of $H_2(X, \mathbb{Z})$ and $\underline{c}$ of $H_4(X, \mathbb{Z})$, respectively. Let $\underline{x}$ and $\underline{y}$ be the dual bases of $H^2(X, \mathbb{Z})$ and $H^4(X, \mathbb{Z})$, respectively. Finally, let $\underline{y}^*$ be the basis of $H^4(X, \mathbb{Z})$ which is via $\gamma_X$ dual to $\underline{y}$.

We find $\partial W_2 \cong \#_{i=1}^b (S^2 \times S^5)$, and $W_4$ is determined by the ambient isotopy class of a framed link of 3-spheres in $\partial W_2$ with $b'$ components. Let $f_k \colon S^3 \times D^4 \longrightarrow \partial W_2$ be the $k$-th component of that link and $g_k := f_{k|S^3 \times \{0\}}$, $k = 1, ..., b'$. In the notation of Section 3.6, we write $[g_k] \in \pi_3(\partial W_2 \setminus \bigcup_{k \neq j} S_j)$ in the form $(l_i^k, i = 1, ..., b, l_{ij}^k, 1 \leq i < j \leq b; \lambda_{kj}, j \neq k)$, $k = 1, ..., b'$. To see the significance of the $l_i^k$ and $l_{ij}^k$, note that, by Remark 3.4, $W_2 \cup H_k^4 \subset X$ is homotopy equivalent to $\big(\bigvee_{i=1}^b S^2\big) \cup_{g_k} D^4$. The cohomology ring of that complex has been computed in Proposition 3.11, so that the naturality of the cup product implies the following formulae for the cup products in $X$:

$$x_i \cup x_j \;=\; \sum_{k=1}^{b'} l_{ij}^k \cdot y_k^*, \quad i \neq j,$$

$$x_i \cup x_i \;=\; \sum_{k=1}^{b'} l_i^k \cdot y_k^*, \quad i = 1, ..., b.$$

Therefore, the $l_i^k$ and $l_{ij}^k$ are determined by $\delta_X$ and $\gamma_X$ (used for computing $\underline{y}^*$), in fact $l_i^k = \gamma_X(\delta(x_i \otimes x_i) \otimes y_k)$ and $l_{ij}^k = \gamma_X(\delta(x_i \otimes x_j) \otimes y_k)$.

To determine the $\lambda_{ij}$ and the framings, we proceed as follows: Look at the embedding $\#_{i=1}^{b}(S^2 \times S^5) \hookrightarrow X$. There exist $b$ embedded 2-spheres $S_1^2, ..., S_b^2$ which represent the basis $\underline{b}$ and which do not meet the given link. Finally, $\#_{i=1}^{b}(S^2 \times S^5)$ obviously possesses a regular neighborhood in $X$ which is homeomorphic to $\#_{i=1}^{b}(S^2 \times S^5) \times D^1$. Thus, we can perform "surgery in pairs" as described in Section 3.1. The result is a 3-connected manifold $X^*$ containing $S^7$. It is by construction the manifold obtained from the framed link in $S^7$ derived from the given one in $\#_{i=1}^{b}(S^2 \times S^5)$ (cf. Section 4.1). We will be finished, once we are able to compare the invariants of $X$ to those of $X^*$. To do so, we look at the *trace of the surgery*, i.e., at $Y = (X \times I) \cup H_1^5 \cup \cdots \cup H_{b'}^5$, the 5-handles being attached along tubular neighborhoods of the $S_i \times \{1\}$ in $X \times \{1\}$. Then, $\partial Y = X \sqcup \overline{X}^*$. The Mayer-Vietoris sequence provides the isomorphisms

$$H_4(X, \mathbb{Z}) \cong H_4\left(X \setminus \bigsqcup_{i=1}^{b'}(S_i \times D^6), \mathbb{Z}\right) \cong H_4(X^*, \mathbb{Z}).$$

Set $H := H_4\left(X \setminus \bigsqcup_{i=1}^{b'}(S_i \times D^6), \mathbb{Z}\right)$. By Lefschetz duality ([5], (28.18)), there is for each $q \in \mathbb{N}$ a diagram (omitting $\mathbb{Z}$-coefficients)

(4)
$$\begin{array}{ccccccc}
H^{q-1}(Y) & \longrightarrow & H^{q-1}(\partial Y) & \longrightarrow & H^q(Y, \partial Y) & \longrightarrow & H^q(Y) \\
\downarrow \cong & & \downarrow \cong & & \downarrow \cong & & \downarrow \cong \\
H_{10-q}(Y, \partial Y) & \longrightarrow & H_{9-q}(\partial Y) & \longrightarrow & H_{9-q}(Y) & \longrightarrow & H_{9-q}(Y, \partial Y)
\end{array}$$

where the left square commutes up to the sign $(-1)^{q-1}$ and the other two commute. We first use it in the case $q = 5$. Look at the commutative diagram

$$\begin{array}{ccc}
H & \xrightarrow{\cong} & H_4(X^*, \mathbb{Z}) \\
\downarrow \cong & & \downarrow \\
H_4(X, \mathbb{Z}) & \longrightarrow & H_4(Y, \mathbb{Z}),
\end{array}$$

in which all arrows are injective, because $H_5(Y, X; \mathbb{Z}) = 0 = H_5(Y, X^*; \mathbb{Z})$ (cf. [17], p. 198). Using the identification $H_4(\partial Y, \mathbb{Z}) = H \oplus H$, we find

(5) $$\operatorname{Im}(H_5(Y, \partial Y; \mathbb{Z})) = \{(y, -y) \in H \oplus H\}.$$

Similar considerations apply to the case $q = 9$. Taking into account that $X^*$ sits in $Y$ with the reversed orientation, (4) shows that the forms $\gamma_X$ and $\gamma_{X^*}$, both defined with respect to the preferred bases, coincide. In the same manner, the pullbacks of $p_1(Y)$ to $H^4(X, \mathbb{Z})$ and $H^4(X^*, \mathbb{Z})$, respectively, agree. Since $X$ and $X^*$ are the boundary components of $Y$, these pullbacks are $p_1(X)$ and $p_1(X^*)$, respectively, and we are done. □

4.3. **Manifolds with given invariants.** One might speculate, especially in view of the classification of E-manifolds in dimension 4 and 6, that the invariants $\delta_X$, $\gamma_X$, and $p_1(X)$ might suffice to classify E-manifolds with $w_2(X) = 0$ in dimension 8. However, Lemma 3.6 shows that these invariants determine only $W_4$ and we still have the choice of an isomorphism in gluing $\#_{i=1}^{b}(S^2 \times S^5)$ to $W_4$, and different gluings may lead to different results. The following example which was communicated to me by J.-C. Hausmann illustrates this phenomenon.

*Example* 4.4. One has $\pi_5(SO(3)) \cong \mathbb{Z}_2$ [32]. Therefore, there are two different $S^2$-bundles over $S^6$, call them $X := S^6 \times S^2$ and $X' := S^6 \widetilde{\times} S^2$. Obviously, $X$ and $X'$ are spin-manifolds with trivial invariants, but one computes $\pi_5(X) \cong \mathbb{Z}_2$ and $\pi_5(X') = \{0\}$.

Fix $b, b'$, and a system $Z$ of invariants in the image of the map $Z^{\text{PL}(\mathscr{C}^\infty)}(b,b')$. As we have seen, $Z$ determines a certain manifold $W_4$ the boundary of which is diffeomorphic to $\#_{i=1}^b (S^2 \times S^5)$ together with a basis $\underline{b}$ for $H_2(\partial W_4, \mathbb{Z})$. Let $\underline{b}_0$ be the natural basis for $H_2(\#_{i=1}^b (S^2 \times S^5), \mathbb{Z})$, and denote by $\text{Iso}_0^{\text{PL}(\mathscr{C}^\infty)}$ the set of piecewise linear (smooth) isomorphisms $f: \#_{i=1}^b (S^2 \times S^5) \longrightarrow \partial W_4$ with $f_*(\underline{b}_0) = \underline{b}$. Our results show that every based piecewise linear (smooth) manifold $(X, \underline{x}, \underline{y})$ with system of invariants $Z$ is piecewise linearly (smoothly) isomorphic to a manifold of the form

$$X(f) := \partial W_4 \cup_f \#_{i=1}^b (S^2 \times S^5) \quad \text{for some } f \in \text{Iso}_0^{\text{PL}(\mathscr{C}^\infty)}$$

with its given bases for $H^2(X(f), \mathbb{Z})$ and $H^4(X(f), \mathbb{Z})$. Conversely, every manifold of the form $X(f)$ is a piecewise linear (smooth) based E-manifold with invariants $Z$.

Now, suppose we are given $f, f' \in \text{Iso}_0^{\text{PL}(\mathscr{C}^\infty)}$, such that $X(f)$ and $X(f')$ are isomorphic as piecewise linear (smooth) based manifolds. We claim that we can find an isomorphism $\varphi: X(f) \longrightarrow X(f')$ with $\varphi(W_4) = W_4$. For this, look at the handle decomposition $W_0 \subset W_2 \subset W_4$. Since $W_0$ is just an embedded 8-disc in $X(f)$ and $X(f')$, respectively, we can choose $\varphi$ with $\varphi(W_0) = W_0$. Let $l \subset \partial W_0$ be the framed link for attaching the 2-handles. Then, $\varphi(l)$ and $l$ are isotopic. Therefore, we can find a level preserving diffeomorphism $\widetilde{\psi}: [-1,1] \times \partial W_0 \longrightarrow [-1,1] \times \partial W_0$ with $\widetilde{\psi}_{|\{\pm 1\} \times \partial W_0} = \text{id}_{\partial W_0}$ and $\widetilde{\psi}_{|\{0\} \times \partial W_0}(\varphi(l)) = l$. If we choose a tubular neighborhood ($\cong [-1,1] \times \partial W_0$) of $\partial W_0$ in $X(f')$, we can use $\widetilde{\psi}$ to define an automorphism $\psi: X(f') \longrightarrow X(f')$ with $\psi(\varphi(l)) = l$. Thus, $\psi \circ \varphi$ maps $W_2$ onto $W_2$. A similar argument shows that we can achieve $\varphi(W_4) = W_4$.

Let $\text{Aut}_0^{\text{PL}(\mathscr{C}^\infty)}(\#_{i=1}^b (S^2 \times D^6))$ be the group of piecewise linear (smooth) automorphisms $g$ of $\#_{i=1}^b (S^2 \times D^6)$ with $H^2(g, \mathbb{Z}) = \text{id}$ and similarly define $\text{Aut}_0^{\text{PL}(\mathscr{C}^\infty)}(W_4)$. Then, we have just established

**Proposition 4.5.** *The set of isomorphy classes of based piecewise linear (smooth) E-manifolds with invariants $Z$ is in bijection to the set of equivalence classes in $\text{Iso}_0^{\text{PL}(\mathscr{C}^\infty)}$ w.r.t. the equivalence relation coming from the group action*

$$\text{Aut}_0^{\text{PL}(\mathscr{C}^\infty)}(W_4) \times \text{Aut}_0^{\text{PL}(\mathscr{C}^\infty)}(\#_{i=1}^b (S^2 \times D^6)) \times \text{Iso}_0^{\text{PL}(\mathscr{C}^\infty)} \longrightarrow \text{Iso}_0^{\text{PL}(\mathscr{C}^\infty)}$$
$$(h, g, f) \longmapsto h_{|\partial W_4} \circ f \circ g^{-1}_{|\#_{i=1}^b (S^2 \times S^5)}.$$

We shall see in Lemma 5.1 that $\text{Aut}_0^{\text{PL}(\mathscr{C}^\infty)}(\#_{i=1}^b (S^2 \times D^6))$ contains the commutator subgroup of $\text{Aut}_0^{\text{PL}(\mathscr{C}^\infty)}(\#_{i=1}^b (S^2 \times S^5))$.

**Corollary 4.6.** *The set of isomorphy classes of based piecewise linear (smooth) E-manifolds with $b_2 = b$ and $b_4 = 0$ is in bijection to the abelian group*

$$\text{Aut}_0^{\text{PL}(\mathscr{C}^\infty)}(\#_{i=1}^b (S^2 \times S^5)) / \text{Aut}_0^{\text{PL}(\mathscr{C}^\infty)}(\#_{i=1}^b (S^2 \times D^6)).$$

I have been informed by experts that the structure of the groups $\mathrm{Aut}_0^{\mathrm{PL}(\mathscr{C}^\infty)}(\#_{i=1}^b(S^2 \times S^5))$ and $\mathrm{Aut}_0^{\mathrm{PL}(\mathscr{C}^\infty)}(\#_{i=1}^b(S^2 \times D^6))$ has not yet been determined and that this would be a rather difficult task. Therefore, we choose the viewpoint of framed links in order to finish our considerations. In Theorem 5.2, we will then use this viewpoint to compute the group $\mathrm{Aut}_0^{\mathrm{PL}}(\#_{i=1}^b(S^2 \times S^5))/\mathrm{Aut}_0^{\mathrm{PL}}(\#_{i=1}^b(S^2 \times D^6))$.

As above, let $(X,\underline{x},\underline{y})$ be a based piecewise linear (smooth) E-manifold with zero second Stiefel-Whitney class and system of invariants $Z_{(X,\underline{x},\underline{y})} = (\delta, \gamma, p)$. We have seen that we can find a framed link $l_X$ of 2-spheres in $X$ which represents the basis $\underline{x}$ and perform surgery along this link in order to get a 3-connected piecewise linear (smooth) based manifold $(X^*,\underline{y})$ together with a framed link $l_{X'}^*$ of 5-spheres in it. If $(X',\underline{x}',\underline{y}',l_{X'})$ is another such object where $(X',\underline{x}',\underline{y}')$ is isomorphic to $(X,\underline{x},\underline{y})$, then we can clearly find an isomorphism $\varphi\colon (X,\underline{x},\underline{y}) \longrightarrow (X',\underline{x}',\underline{y}')$ with $\varphi(l_X) = l_{X'}$. Such an isomorphism $\varphi$ yields, after surgery, an isomorphism $\varphi^*\colon (X^*,\underline{y}) \longrightarrow (X'^*,\underline{y}')$ with $\varphi^*(l_{X'}^*) = l_{X'^*}^*$. In particular, the manifold $(X^*,\underline{y})$ is determined up to piecewise linear (smooth) isomorphy. We call it the *type of* $(X,\underline{x},\underline{y})$. Note that this notion matters only in the smooth case, by Theorem 2.2. To summarize, we note

**Proposition 4.7.** *The set of isomorphy classes of based piecewise linear (smooth) E-manifolds of type $(X^*,\underline{y})$ is in bijection to the set of equivalence classes of framed links of 5-spheres in $X^*$ where two such links $l$ and $l'$ are considered* equivalent, *if there is a piecewise linear (smooth) automorphism $\varphi^*\colon (X^*,\underline{y}) \longrightarrow (X^*,\underline{y})$ with $\varphi^*(l) = l'$.*

*Example* 4.8. The group $\mathbb{Z}_2^{\oplus b}$ acts freely on the set of isotopy classes of framed links of $b$ spheres of dimension 5 in $X^*$ by altering the framings of the components. Note that the two possible framings of the trivial bundle on a 5-sphere are distinguished by the fact that one extends over $D^6$ and the other does not. This property is preserved under piecewise linear homeomorphisms, so that we conclude that $\mathbb{Z}_b^{\oplus b}$ acts also freely on the set of equivalence classes of framed links of $b$ spheres of dimension 5 in $X^*$.

Note that this completes the classification of Spin-E-manifolds of dimension eight with second Betti number one.

Let us look at manifolds of type $S^8$. We claim that two framed links $l$ and $l'$ of 5-spheres are equivalent in the above sense, if and only if they are isotopic. Clearly, after replacing $l$ and $l'$ by isotopic links, we may assume that both of them are contained in the Southern hemisphere and that $\varphi^*$ is the identity on the Northern hemisphere. Now, choose a representative $\varphi^\dagger$ for the isotopy class of $\varphi^{*-1}$ which is the identity on the Southern hemisphere. Then, $\varphi^\dagger \circ \varphi^*$ is isotopic to the identity and carries $l$ into $l'$.

For differentiable manifolds, the operation $X \longmapsto X \# \Sigma$, $\Sigma$ an exotic 8-sphere, establishes a bijection between the set of isomorphy classes of based smooth E-manifolds of type $S^8$ and the set of isomorphy classes of based smooth E-manifolds of type $\Sigma$. We conclude

**Corollary 4.9.** i) *The set of isomorphy classes of based piecewise linear E-manifolds with $b_2 = b$ and $b_4 = 0$ is in bijection to the group $\mathrm{FL}_b = \mathrm{L}_b \oplus \bigoplus_{i=1}^b \mathbb{Z}_2$.*

ii) *The set of isomorphy classes of based smooth E-manifolds with $b_2 = b$ and $b_4 = 0$ is in bijection to the group $\vartheta^8 \oplus \mathrm{FL}_b$.*

Finally, we have to deal with those manifolds for which the cup form $\delta$ is trivial. Our investigations in Section 3.6 and 4.2 show that the framed link of 3-spheres in $\partial W_2$ can be chosen to be contained in a small disc. In other words, a manifold $X$ with $\delta_X \equiv 0$ is piecewise linearly (smoothly) isomorphic $X^\dagger \# X^*$ where $X^*$ is the type of $X$ and $b_4(X^\dagger) = 0$. As our surgery arguments above reveal, an isomorphism between $X^\dagger \# X^*$ and $X'^\dagger \# X'^*$ can be chosen of the form $\varphi^\dagger \# \varphi^*$ where $\varphi^\dagger \colon X^\dagger \longrightarrow X'^\dagger$ and $\varphi^* \colon X^* \longrightarrow X'^*$ are isomorphisms. Therefore, the set of isomorphy classes of based piecewise linear E-manifolds of type $X^*$ with $b_2 = b$ is in bijection to the set of isomorphy classes of based piecewise linear E-manifolds with $b_2 = b$ and $b_4 = 0$. The same goes for differentiable manifolds of type $X^*$, if $X^*$ is not diffeomorphic to $X^* \# \Sigma$, $\Sigma$ an exotic 8-sphere. Otherwise, we have to divide by the action of $\vartheta^8$. This observation together with Corollary 4.9 settles Theorem 2.4. □

## 5. The structure of the group $\mathrm{Aut}_0^{\mathrm{PL}}\big(\#_{i=1}^b(S^2 \times S^5)\big) / \mathrm{Aut}_0^{\mathrm{PL}}\big(\#_{i=1}^b(S^2 \times D^6)\big)$

The aim of this section is to prove that, first, $\mathrm{Aut}_0^{\mathrm{PL}}\big(\#_{i=1}^b(S^2 \times S^5)\big) / \mathrm{Aut}_0^{\mathrm{PL}}\big(\#_{i=1}^b(S^2 \times D^6)\big)$ is an abelian group, and, second, that it is, in fact, isomorphic to the group $\mathrm{FL}_b$ defined before. This result should be of some independent interest, especially because the group $\mathrm{FL}_b$ is by Haefliger's work quite well understood. For $b = 1$, we refer to [20] for more specific information.

We begin with the elementary

**Lemma 5.1.** *Let $k \in \mathrm{Aut}_0^{\mathrm{PL}}\big(\#_{i=1}^b(S^2 \times S^5)\big)$ be a commutator. Then, $k$ extends to an automorphism of $\#_{i=1}^b(S^2 \times D^6)$.*

*Proof.* For the proof, we depict $\#_{i=1}^b(S^2 \times S^5)$ as follows: Let $V_i$, $i = 1, ..., b$, be $b$ copies of $S^2 \times D^6$, and we join $V_i$ and $V_{i+1}$ by a tube $T_i \cong [-1, 1] \times D^7$, $i = 1, ..., b-1$. The result is a manifold $W$ whose boundary is isomorphic to $\#_{i=1}^b(S^2 \times S^5)$. We make the following normalizations: Write $\partial V_i$ as $(S^2 \times D_+^i) \cup (S^2 \times D_-^i)$, let $n_i$ and $s_i$ be the centers of $D_+^i$ and $D_-^i$, respectively, and set $S_+^i := S^2 \times n_i$ and $S_-^i := S^2 \times s_i$, $i = 1, ..., b$. Choose furthermore points $e_i \neq w_i$ in $(S^2 \times D_+^i) \cap (S^2 \times D_-^i)$, $i = 1, ..., b$, and suppose that $\{-1\} \times D^7 \subset T_i$ is attached to a disc around $w_i$ in $\partial V_i$ and $\{1\} \times D^7 \subset T_i$ to a disc around $e_{i+1}$ in $\partial V_{i+1}$, $i = 1, ..., b-1$. Set $T := \bigsqcup_{i=1}^{b-1} T_i$.

Now, let $k = f \circ g \circ f^{-1} \circ g^{-1}$ with $f, g \in \mathrm{Aut}_0^{\mathrm{PL}}\big(\#_{i=1}^b(S^2 \times S^5)\big)$. As $H_2(h, \mathbb{Z})$ is the identity for every element $h \in \mathrm{Aut}_0^{\mathrm{PL}}\big(\#_{i=1}^b(S^2 \times S^5)\big)$ and $S_\pm^i$, $i = 1, ..., b$, both represent the same basis for $H_2(\partial W, \mathbb{Z})$, $h$ is isotopic to a map $h'$ which satisfies either assumption (A) or (B) below.

- (A): $h'$ is trivial on a tubular neighborhood of $S_+^i$ which contains $(S^2 \times D_+^i) \setminus \mathrm{Int}(T)$, $i = 1, ..., b$.
- (B): $h'$ is trivial on a tubular neighborhood of $S_-^i$ which contains $(S^2 \times D_-^i) \setminus \mathrm{Int}(T)$, $i = 1, ..., b$.

Next, replace $f$ by an isotopic map $f'$, satisfying (A), and $g$ by an isotopic map $g'$, satisfying (B). Then, $k'$ is isotopic to $f' \circ g' \circ f'^{-1} \circ g'^{-1}$. The map $k'$ is the identity outside $\mathrm{Int}(\partial T)$. It is, furthermore, the identity on a collar of $(\{-1\} \sqcup \{1\}) \times S^6$ in $R_i := [-1, 1] \times S^6 \subset \partial T_i$, $i = 1, ..., b-1$. Let $k_i'$ be the restriction of $k'$ to $R_i$, $i = 1, ..., b$. We know that each $k_i'$ is the identity on a collar of $\{-1, 1\} \times S^6$ in $R_i$. Thus, we extend every $k_i'$ to a PL automorphism $\widetilde{k}_i$ of $D^7 \times \{-1\} \cup R_i \cup D^7 \times \{1\} \cong S^7$ through $\mathrm{id}_{D^7 \times \{-1\} \cup D^7 \times \{1\}}$. Now, by [27], Lemma 1.10, p. 8, $\widetilde{k}_i$ extends to an automorphism $\kappa_i$ of $D^8 \cong D^7 \times [-1,1]$, $i = 1,...,b$. Thus, the maps $\mathrm{id}_{V_i}$ and $\alpha_i$, $i = 1,...b$, glue to an automorphism of $\#_{i=1}^b(S^2 \times D^6)$ the restriction of which to the boundary is just $k'$. □

This Lemma shows that the group $\mathrm{Aut}_0^{\mathrm{PL}}\bigl(\#_{i=1}^b(S^2 \times D^6)\bigr)$ is normal in the group $\mathrm{Aut}_0^{\mathrm{PL}}\bigl(\#_{i=1}^b(S^2 \times S^5)\bigr)$, and that $\mathrm{Aut}_0^{\mathrm{PL}}\bigl(\#_{i=1}^b(S^2 \times S^5)\bigr)/\mathrm{Aut}_0^{\mathrm{PL}}\bigl(\#_{i=1}^b(S^2 \times D^6)\bigr)$ is abelian. Moreover, in Section 4.3, we have already defined a set theoretic bijection

$$\beta\colon \mathrm{Aut}_0^{\mathrm{PL}}\bigl(\#_{i=1}^b(S^2 \times S^5)\bigr)/\mathrm{Aut}_0^{\mathrm{PL}}\bigl(\#_{i=1}^b(S^2 \times D^6)\bigr) \longrightarrow \mathrm{FL}_b.$$

**Theorem 5.2.** *The map $\beta$ is a group isomorphism.*

*Proof.* Since $\beta$ is bijective, we have to verify that $\beta$ is a homomorphism. In order to do so, we will construct a group $\mathbb{G}$ together with surjective homomorphisms

$$\chi_1\colon \mathbb{G} \longrightarrow \mathrm{Aut}_0^{\mathrm{PL}}\bigl(\#_{i=1}^b(S^2 \times S^5)\bigr)/\mathrm{Aut}_0^{\mathrm{PL}}\bigl(\#_{i=1}^b(S^2 \times D^6)\bigr)$$

and

$$\chi_2\colon \mathbb{G} \longrightarrow \mathrm{FL}_b,$$

such that $\chi_2 = \beta \circ \chi_1$. This will clearly settle the claim.

Before we define $\mathbb{G}$, we recall some constructions and conventions from [11]. Let $S^8 = \{(x_0,...,x_9) \in \mathbb{R}^9 \,|\, x_0^2 + \cdots x_9^2 = 1\}$ be the unit sphere, write $S^8 = D_+^8 \cup D_-^8$, and let $\sigma\colon S^8 \longrightarrow S^8$ be the reflection at $S^7 = D_+^8 \cap D_-^8$, interchanging the Northern and the Southern hemissphere. First, let $S_b := (S_1^5,...,S_b^5)$ be a "standard link" in $S^8$ defined as follows: Fix real numbers $-1/2 < a_1 < \cdots < a_b < 1/2$, and set

$$S_i^5 \;:=\; \bigl\{(x_0,...,x_9) \in S^8 \,|\, x_6 = x_7 = x_8 = 0,\ x_9 = a_i\bigr\}.$$

We choose, furthermore, framings $\tau_i\colon S_i^5 \times D^3 \longrightarrow S^8$ which extend over $D^6$, such that $\tau_i(D_{i,\pm}^5 \times D^3) \subset D_\pm^8$ and $\sigma \circ \tau_i = \tau_i \circ (\sigma_{|S_i^5} \times \mathrm{id}_{D^3})$, $i = 1,...,b$. Let $l_b^0$ be the resulting framed link in $S^8$ with $l_{b,\pm}^0 := l_b^0 \cap D_\pm^8$. Recall from Section 1 of [11] that

1. Every framed link $l$ of $b$ five dimensional spheres in $S^8$ is isotopic to a link $l'$, such that either (A) $l' \cap D_+^8 = l_{b,+}^0$ or (B) $l' \cap D_-^8 = l_{b,-}^0$.
2. If $l_1$ satisfies (A) and $l_2$ satisfies (B), then $l_1 + l_2$ is represented by the link $l$ with $l \cap D_+^8 = l_2 \cap D_+^8$ and $l \cap D_-^8 = l_1 \cap D_-^8$.

Note that, if we perform surgery along $l_b^0$, we get a manifold $W = W_+ \cup W_-$ which is isomorphic to $\#_{i=1}^b(S^2 \times S^6)$, and $W_\pm := \bigl(D_\pm^8 \setminus \mathrm{Int}(l_b^0)\bigr) \cup \bigl(\bigsqcup_{i=1}^b(S_i^2 \times D_\pm^6)\bigr)$ is canonically isomorphic to $\#_{i=1}^b(S^2 \times D^6)$. For the rest of the proof, we will use the description of $\#_{i=1}^b(S^2 \times S^5)$ as $\partial W_+ = \partial W_-$. Set

$$\mathbb{G} \;:=\; \Bigl\{\, \mathrm{PL\text{-}maps}\ f\colon S^7 \setminus \mathrm{Int}(l_b^0) \longrightarrow S^7 \setminus \mathrm{Int}(l_b^0) \,|\, f_{|\mathrm{boundary}} = \mathrm{id}\Bigr\}.$$

For every $f \in \mathbb{G}$, we define $\varphi(f)\colon \#_{i=1}^b(S^2 \times S^5) \longrightarrow \#_{i=1}^b(S^2 \times S^5)$, by extending $f$ through the identity on $\bigsqcup_{i=1}^b(S_i^2 \times D^5)$. Similarly, define $\psi(f)\colon S^7 \longrightarrow S^7$. Obviously,

$$\begin{array}{rccc}\chi_1\colon & \mathbb{G} & \longrightarrow & \mathrm{Aut}_0^{\mathrm{PL}}\bigl(\#_{i=1}^b(S^2 \times S^5)\bigr)/\mathrm{Aut}_0^{\mathrm{PL}}\bigl(\#_{i=1}^b(S^2 \times D^6)\bigr) \\ & f & \longmapsto & [\varphi(f)]\end{array}$$

is a surjective homomorphism.

Next, we associate to $f \in \mathbb{G}$ an element $\chi_2(f) \in \mathrm{FL}_b$ as follows: First, we define $\Sigma(f) := D_+^8 \cup_{\psi(f)} D_-^8$ and the link $l'(f) := l_{b,+}^0 \cup_{\psi(f)} l_{b,-}^0$. Then, we choose a piecewise linear homeomorphism $F\colon \Sigma(f) \longrightarrow S^8$ and set $l_F(f) := F(l'(f))$. We have checked before that the isotopy class of $l_F(f)$ does not depend on the chosen homeomorphism, so that $\chi_2(f) := [l_F(f)] \in \mathrm{FL}_b$ is well defined. To see that $\chi_2 \colon \mathbb{G} \longrightarrow \mathrm{FL}_b$ is a homomorphism, let $f, f'$ be in $\mathbb{G}$. Choose extensions $\overline{\psi} \colon D_+^8 \longrightarrow D_+^8$ and $\overline{\psi}' \colon D_-^8 \longrightarrow D_-^8$ of $\psi(f)$ and $\psi(f')$, respectively. We then define $F \colon \Sigma(f) \longrightarrow S^8$ as $\overline{\psi}$ on $D_+^8$ and as the identity on $D_-^8$, $F' \colon \Sigma(f) \longrightarrow S^8$ as the identity on $D_+^8$ and $(\overline{\psi}')^{-1}$ on $D_-^8$, and $F'' \colon \Sigma(f' \circ f) \longrightarrow S^8$ as $\overline{\psi}$ on $D_+^8$ and $(\overline{\psi}')^{-1}$ on $D_-^8$. Then, the link $l_F(f)$ satisfies (B), the link $l_{F'}(f')$ satisfies (A), and (2) above shows that $[l_{F''}(f' \circ f)] = [l_{F'}(f')] + [l_F(f)]$.

Finally, for given $f \in \mathbb{G}$, we can perform surgery on $\Sigma(f)$ along $l'(f)$. The result is $W_+ \cup_{\varphi(f)} W_-$. Reading this backwards means nothing else but $\beta(\chi_1(f)) = \chi_2(f)$ and we are done. □

UNIVERSITÄT GH ESSEN, FB6 MATHEMATIK & INFORMATIK, D-45141 ESSEN, GERMANY
*E-mail address*: `alexander.schmitt@uni-essen.de`